\documentclass[12pt]{article}
\usepackage{amsmath}
\usepackage{amssymb}
\include{amssymb}

\newtheorem{theorem}{Theorem}[section]

\newtheorem{definition}{Definition}[section]

\newtheorem{proposition}{Proposition}[section]
\newtheorem{remark}{Remark}[section]
\newtheorem{example}[theorem]{Example}

\def\refs{\smallskip\hangindent=25pt\hangafter=1}

\begin{document}
\pagestyle{plain}
\pagenumbering{arabic}

\title{\bf ON CHAPTER XII IN CARTAN'S "LE\c{C}ONS SUR LA
G\'{E}OM\'{E}TRIE DES ESPACES DE RIEMANN". }
\author{Vic Patrangenaru \thanks{Research
supported by NSF-DMS-0805977 and NSA-MSP-H98230-08-1-0058}\\
Florida State University\\
vic@stat.fsu.edu }
\date{ }
\maketitle

\begin{abstract} One shows that Cartan's method of adapted frames in Chapter XII of his famous treatise of Riemannian geometry,
leads to a classification theorem of homogeneous Riemannian
manifolds. Examples of classification in 3D dimensions obtained by
Cartan are given using this powerful method.
\end{abstract}
\refs AMS 1991 {\it subject classifications.} Primary 53C30

\section{Introduction.}
The first edition of Cartan's Riemannian geometry textbook ( some
think, the most important text on Riemannian geometry to this day),
was based on Cartan's lectures in 1925. An augmented edition
appeared in 1946 [second edition, Gauthier-Villars, Paris, 1946]. In
view of the great influence which this work has exerted on the
subsequent development of Differential Geometry. The book was
translated in English and in Russian, and carefully reviewed in Math
Reviews in 1985. That review, which was presenting Cartan's work to
a worldwide audience, had an important omission: Cartan's
contribution in Chapter XII: Groups of Isometries of a Riemannian
Manifold. To these days the elegance and the richness in ideas {\it
adapted frames} that were at the center of that book chapter, are
still little known to the general community of mathematicians, and
surprisingly enough, are rarely quoted by Riemannian geometry
experts. In this paper we will go to over that Chapter with having
also on mind more recent developments of Cartan's ideas.
\section{ Simply transitive groups of isometries of
Riemannian manifolds}

The objective of Chapter XII in \cite{C} was to (i) {\it determine
those Lie groups that act as groups of isometries of a Riemannian
manifold}, and (ii) {\it for such a Lie group, find the Riemannian
manifolds
admitting that group as a group of isometries.} \\
Assume $\bf K$ is a Lie group, subgroup of the group $\bf I_gM$ of
isometries of the Riemannian manifold $\bf (M, g).$ For each point
$x \in  \bf M,$ the {\it orbit } of $x$ ( called by Cartan {\it
trajectory}) is $\bf K(x) = \{k(x), k \in \bf K.\}$ Cartan
introduces the method of {\it adapted frames}. For each $u \in \bf
OM_x,$ Cartan defines the {\it adapted frame field} to be frame
field defined on the orbit $\bf K(x)$ as follows : \begin{equation}
\label{adapted} U(kx) = (d_xk)(u), \forall k \in \bf K.
\end{equation} The group $\bf K$ is {\it transitive } if $\bf K(x) =
M,$ which is same as saying that $(\bf M, g)$ is a {\it Riemannian
homogeneous space.} Cartan considers first the particular case when
the group $\bf K$ is {\it simply transitive}. In this case in the
{\it isotropy group } $\bf G = K_x,$ subgroup of elements in $\bf K$
that keep the point $x$ fixed, is trivial, and the manifold $\bf M$
and $\bf K$ are in a one to one correspondence, thus $\bf M$ has a
Lie group structure. Moreover the adapted frame field $U$ defined in
\eqref{adapted} determines the Riemannian metric on $\bf M,$ since
the if we define $\Theta_u$ to be the {\it dual form} of the
orthoframe field $U,$ then the Riemannian metric $g$ is given by
\begin{equation}\label{g}
g = \|\Theta_u\|^2.
\end{equation}.
Assume $\bf M$ has dimension $n$ and let
\begin{equation}\label{simply} \Theta_u (k(x) = \theta^1(k(x)) U_1(k(x)) +\dots
+ \theta^n(k(x)) U_n(k(x)).\end{equation} Then since the form
$\Theta_u$ is invariant ( $\forall k \in \bf K, L_k^*\Theta_u =
\Theta_u$ ), it follows that $d\Theta_u$ is also invariant,
therefore \begin{equation} \label{constants} d\theta^i = {1\over 2}
C_{jk}^i\theta^j\wedge \theta^k, \end{equation} where $C_{jk}^i$ are
constants with $C_{jk}^i + C_{kj}^i = 0.$ \\
On the other hand the {\it connection forms} associated with the
frame filed $\Theta_u$ are defined by \begin{equation}
\label{connection} d\theta^i + \theta^i_j\wedge\theta^j =0,
\theta^i_j + \theta^j_i = 0, \forall i, j =1,\dots, n.
\end{equation}
From \eqref{constants}, \eqref{connection} Cartan obtains
\begin{equation}\label{christoffel}
\theta^i_j = {1\over 2} C_{jk}^i\theta^k. \end{equation} On the
other hand the {\it curvature forms} associated with the frame filed
$\Theta_u$ are defined by \begin{equation} \label{curvature}
\Omega^i_j = d\theta^i_j + {1\over 2} \theta^i_k\wedge\theta^k_j =0,
\Omega^i_j + \Omega^j_i = 0, \forall i, j =1,\dots, n,
\end{equation}
and the {\it Riemann-Christoffel} curvature coefficients $R^i_{jkr}$
with respect to this frame are given by
\begin{equation} \label{curvature} \Omega^i_j ={1\over 2}
R^i_{jkr} \theta^k\wedge\theta^r, R^i_{jkr} + R^i_{jrk} = 0, \forall
i, j, k, r =1,\dots, n,
\end{equation}
and since the we work in an orthogonal frame filed, $R^i_{jkr} =
R_{ijkr}$ From equations \eqref{connection}-\eqref{curvature} Cartan
obtains
\begin{equation}\label{curvature-Lie}
R_{ijkr} = -{1\over 4}C^m_{ij}C_{krm}.
\end{equation}

\begin{example} In Chapter XII Cartan considers the particular example of the unit sphere $S^3=\{x = (x^1, x^2, x^3, x^4) \in \mathbb R^4, \|x\|^2
=1,$ and considers the Pfaff forms
\begin{eqnarray}\label{sigma}
\sigma_1=-x^2dx^1+x^1dx^2-x^4dx^3+x^3dx^4 \nonumber \\
\sigma_2=-x^3dx^1+x^4dx^2+x^1dx^3-x^2dx^4 \\
\sigma_3=-x^4dx^1-x^3dx^2+x^2dx^3+x^1dx^4  \nonumber
\end{eqnarray}
that correspond to an adample frame field corresponding to a simply
transitive group of isometries. Milnor (\cite{M}) showed that any
left invariant matric on $S^3$ is of the form $g_{\lambda},$ where
 \begin{equation}
 g_{\lambda} = 4((\lambda_2\lambda_3)^{-1}\sigma_1^2+
(\lambda_1\lambda_3)^{-1}\sigma_2^2+
(\lambda_1\lambda_2)^{-1}\sigma_3^2).
\end{equation}
\end{example}
\section{Riemannian manifolds that admit a multiply transitive group of isometries.}
The second part of Chapter XII in \cite{C} is dedicated to
transitive groups of isometries $\bf K$ of a Riemannian manifold
$(\bf M, g)$ for which the dimension of the isotropy group $\bf H =
K_x$ is positive. As a starting point, Cartan  the {\it isotropic
representation} $\mathfrak g = \lambda_u(\mathfrak h)$ (\cite{KN1})
of the {\it isotropy algebra $\mathfrak h$} as a Lie subalgebra of
$\mathfrak{o}(n).$ Note that the Lie group $\bf K$ is embedded in
the orthoframe bundle $\bf O_gM,$ via
\begin{equation}\label{embedding}
\phi_u(k) = (d_xk)(u), \forall k \in \bf K.
\end{equation}
Let $\theta \in \Omega^1(\bf O_gM,\mathbb R^n), \omega \in
\Omega^1(\bf O_gM,\mathfrak{o}(n))$ be the dual, respectively the
Levi-Civita connection form of $(M,g).$ The vector valued
differential forms $\theta_u = \phi_u^*\theta, \omega_u =
\phi_u^*\omega$ are left-invariant differential forms on $\bf K,$
with $\theta_u \in \Omega^1(\bf K,\mathbb R^n), \omega_u \in
\Omega^1(\bf K,\mathfrak{o}(n)),$ and if $\mathfrak{g}^{\perp}$ is
the orthocomplement of $\mathfrak g$ in $\mathfrak{o}(n),$ then the
left invariant forms $\theta_u, \omega_{u, \mathfrak g}$ on $\bf K$
are linearly independent, and the Levi-Civita connection form, then
the $\phi_u$ pull-back of the structure equations,
\begin{equation} \label{structural}
d\theta + \omega \wedge \theta = 0, \Omega = d\omega +
\frac{1}{2}\omega \wedge \omega,
\end{equation}
yield the Maurer-Cartan equations of $\bf K.$
\begin{equation} \label{Maurer-Cartan}
d\theta_u + \omega_u \wedge \theta_u = 0, \Omega_u = d\omega_u +
\frac{1}{2}\omega_u \wedge \omega_u,
\end{equation}
If $\omega_{u} = \omega_{u, \mathfrak g}\oplus \omega_{u,
\mathfrak{g}^{\perp}}$ is the decomposition of $\omega_{u}$ with
respect to the splitting $\mathfrak{o}(n) = \mathfrak g \oplus
\mathfrak{g}^{\perp},$ the $\mathfrak{g}^{\perp}$-component
$\omega_{u, \mathfrak{g}^{\perp}}$ depends on $\theta_u$ only as
shown by Cartan, therefore there is a $\mathfrak g$-invariant linear
map $\Gamma: \mathbb R^n \to \mathfrak{g}^{\perp},$ such that
$\omega_{u, \mathfrak{g}^{\perp}} = \Gamma (\theta_u).$ Similarly,
if $\Omega_u = \Omega_{u, \mathfrak g}\oplus \Omega_{u,
\mathfrak{g}^{\perp}},$ there is a $\mathfrak g$-invariant bilinear
map $\bar{\Omega}:\mathbb R^n \times \mathbb R^n : \to \mathfrak g,$
such that $\Omega_{u, \mathfrak g}= \bar{\Omega}(\theta_u,
\theta_u).$ The Lie algebra equations, dual to \eqref{Maurer-Cartan}
are:
\begin{eqnarray}\label{cartan-triple}
T(X,Y) = \Gamma(Y)X - \Gamma(X)Y \nonumber \\
\tilde{\Omega}(X,Y) = \overline{\Omega} - [\Gamma(X), \Gamma(Y)]_{\frak{g}} \nonumber \\
\Big[\xi,\eta \Big] = [\xi,\eta] \\
\Big[ \xi, X \Big] = \xi(X) \nonumber \\
\Big[ X, Y \Big] = - T(X,Y) - \tilde{\Omega}(X,Y),  \nonumber
\end{eqnarray}
leading us to the following definition
\begin{definition}
An $n-$dimensional {\it Cartan triple } in a triple $(\mathfrak{g},
\Gamma, \overline{\Omega})$ where $ \Gamma : \mathbb{R}^n
\rightarrow \mathfrak{g}^{\perp}$ is a $\frak{g}$-invariant linear
map and  $$\overline{\Omega} : \mathbb{R}^n \times \mathbb{R}^n
\rightarrow \frak{g}$$ is a $\frak{g}$-invariant bilinear map such
that if we define $T, \tilde{\Omega}$ and  the ``taller" bracket $
\Big[ , \Big]$ by \eqref{cartan-triple}, then $ \Big[ , \Big]$
yields a Lie algebra structure on $\frak{k}(\frak{g},\Gamma,
\overline{\Omega})=\frak{g} \oplus \mathbb{R}^n,$ provided some
identities in $\tilde{\Omega}$ and $\Gamma,$ resulting
from the Jacobi identities for that Lie algebra hold true. \\
Here $[,]$ is the commutator of two matrices and $\xi_\frak{g}$ is
the $\frak{g}$-component of $\xi$ with respect to the decomposition
$\frak{o}(n)=\frak{g}\oplus \frak{g}^{\perp}.$
\end{definition}

\subsection{ Maximal closed Cartan triples and full groups of
isometries of Riemannian homogeneous spaces}

Cartan triples were introduced by Patrangenaru \cite{P1}, who used
them to classify metrically Riemannian homogeneous spaces. Note that
for each Cartan triple $(\frak{g},\Gamma, \overline{\Omega}),
\frak{g}$ is a Lie subalgebra of $\frak{k}=\frak{k}(\frak{g},\Gamma,
\overline{\Omega})$. Let $\bf K$ be the simply connected Lie group
of Lie algebra $\frak{k}$ and let $\bf G$ be the connected Lie
subgroup of $\bf K$, whose Lie algebra is $\frak{g}.$ The Cartan
triple $(\frak{g},\Gamma, \overline{\Omega})$
is said to be {\it closed} if $\bf G$ is a closed subset of $\bf K.$ \\
\indent Assume $\mathcal{R}_1, \mathcal{R}_2$ are two Cartan
triples, such that for $j=1,2, \mathcal{R}_j = (\frak{g}_j,\Gamma_j,
\overline{\Omega}_j).$ We say that $ \mathcal{R}_1 \leq
\mathcal{R}_2$ if there is a vector subspace $\frak{a}$
of $\frak{g}_2$ such that \\

\begin{equation}
\frak{g}_1 \subseteq \frak{g}_2=\frak{g}_1\oplus \frak{a},
\frak{g}_1^{\perp}=\frak{g}_2^{\perp}\oplus \frak{a}
\end{equation}
and with respect to the decompositions in (3) we have
\begin{eqnarray}
\Gamma_1=\Gamma_2\oplus\Gamma_{\frak{a}},\\
\overline{\Omega}_2=\overline{\Omega}_1\oplus\overline{\Omega}_{\frak{a}},
\end{eqnarray}
where $\Gamma_{\frak{a}}$ is the $\frak{a}$-component of $\Gamma_1$
and $\overline{\Omega}_{\frak{a}}$ is the $\frak{a}$-component of
$\overline{\Omega}_2$
with respect to the corresponding decomposition in (3).\\
\indent Let $\mathcal{C}_n$ be the set of all $n$-dimensional Cartan
triples and let $\mathcal{M}_n$ be the set of  maximal Cartan
triples in $(\mathcal{C}_n,\leq).$ The orthogonal group $\bf O(n)$
acts on the right on $\mathcal{C}_n,$ leaving $\mathcal{M}_n$
invariant. This action $A_n$ is given by $A_n((\frak{g},\Gamma,
\overline{\Omega}),a)= (\frak{g}',\Gamma', \overline{\Omega}')$,
where
\begin{eqnarray}
\frak{g}'=Ad(a^{-1})\frak{g},\\
\Gamma'(\cdot)=Ad(a^{-1})\Gamma(a(\cdot)),\\
\overline{\Omega}'(\cdot,\cdot)=Ad(a^{-1})\overline{\Omega}(a(\cdot),a(\cdot)).
\end{eqnarray}
The following result is an immediate consequence of the main result
in \cite{P1}.
\begin{theorem}\label{clssif}(\cite{P1}). There is a one to one correspondence between isometry
classes of simply connected $n$- dimensional Riemannian homogeneous
spaces and $A_n$- orbits of closed Cartan triples in
$\mathcal{M}_n.$
\end{theorem}
In this correspondence, if $(\bf M,g)$ is a Riemannian homogeneous
space, if $u\in \bf O_gM_x$ is a fixed orthoframe at a given point $
x\in \bf M,$ if $\frak{k}(\bf M)_x$ is the Lie algebra of Killing
vector fields on $M$ vanishing at $x$, and if $\lambda_u$ is the
linear isotropic representation, then one may take
$\frak{g}=\lambda_u(\frak{k}(\bf M)_x),$ and $T$ and
$\tilde{\Omega}$ in equation (2) are the {\it torsion and the
$\frak{g}$-component of the curvature} of the {\it Ambrose-Singer
connection}. Conversely a simply connected Riemannian homogeneous
space corresponding to a closed maximal Cartan triple is usually
called a {\it geometric realization} of that Cartan triple. As a
manifold is the quotient space $\bf K/G,$ where $\bf K$ is the
simply connected Lie group of Lie algebra $\frak{k}$ and $G$ is the
Lie subgroup corresponding to $\frak{g}.$ The subgroup $G$ is simply
connected by the exact homotopy sequence of the fibration $\bf G
\subseteq \bf K \rightarrow \bf K/G.$ The metric is defined by the
Euclidean norm in $\mathbb{R}^n,$  and, by the definition of the
Cartan triple is
$ad(\frak{g})$-invariant, hence $\bf K/G$ is a reductive homogeneous space \cite{KN}.\\
\begin{remark}\label{algo} Theorem \ref{clssif} in the line of Cartan's algorithmic program of classification
of Riemannian homogeneous spaces in dimension $n.$ That program was
based on the following three steps: \\(i) Determine all the closed
Lie subgroups of $O(n).$ \\ (ii) For a given Lie subgroup $\bf G
\subseteq O(n),$ determine all the Lie groups $\bf K$ acting as
transitive groups of isometries of a Riemannian manifold with
isotropy group $\bf G,$ and \\
(iii) Given a Lie group $\bf K$ as in (ii), determine all the
homogeneous spaces $\bf(M, g)$ admitting $\bf K$ as a group of
isometries.
\end{remark}
\section{ The 3D Riemannian homogeneous spaces with a multiply transitive group of isometries.}
Let $(e_j^i)_{i,j=1,..,n}$ be the natural basis of
$\frak{gl}(n,\mathbb{R});$ the matrices $e_j^i$ act on the natural
basis of $\mathbb{R}^n$ as linear endomorphisms by:
\begin{equation}
(e_j^i)(e_k)=\delta_k^ie_j , \forall i,j,k=\overline{1,n}.
\end{equation}
An orthogonal basis of $\frak{o}(n)$ with respect to the Killing
form is given by
\begin{equation}
f_j^i=e_j^i-e_i^j, \forall(i,j), 1\leq j <i \leq n .
\end{equation}
Cartan considered in a separate section the case $n=3.$ As a first
step in the Cartan triple method, note that any nontrivial, proper
Lie subalgebra of $\frak{o}(3)$ is conjugated to $\frak{g}_1=
\mathbb{R}f_1^2.$ Therefore, in the multiply transitive case, by
Theorem \ref{clssif} it suffices to list only $A_3$-orbits of
the following ``special" Cartan triples: \\
- $(\frak{o}(3), 0, \overline{\Omega}),$ and \\
- $(\frak{g}_1, \Gamma,  \overline{\Omega}).$ \\
An $(\frak{o}(3), 0, \overline{\Omega})$ Cartan triple corresponds
to a Riemannian manifold whose group of isometries is six
dimensional. In this case a geometric realization of such a Cartan
triple has constant curvature (\cite{K}, Theorem 3.1). A three
dimensional simply connected Riemannian manifold $(M,g)$ with
$dim(I(M,g))=6$ and positive sectional curvature is isometric to a
round sphere $\mathbb{S}^3_R$ of radius $R>0$, the
$\frak{o}(3)$-curvature of which has the form
$\overline{\Omega}=R^{-2}\Omega_1.$ Here
\begin{equation}
\Omega_1(x,y)=(x^1y^2-x^2y^1)f_1^2+(x^1y^3-x^3y^1)f_1^3+(x^2y^3-x^3y^2)f_2^3.
\end{equation}
\'Elie Cartan's application of his method of adapted frames in 3D,
in our terminology amounts to listing the $(\frak{g}_1, \Gamma,
\overline{\Omega})$ - Cartan triples. Cartan's results are as
follows:
\begin{eqnarray}\label{3tripoles}
\Gamma(e_1) = af_1^3 +bf_2^3 , \nonumber \\
\Gamma(e_2) = -bf_1^3 +af_2^3 , \nonumber \\
\Gamma(e_3) = 0, \\
\overline{\Omega}(e_1,e_2)=kf_1^2, \nonumber \\
\overline{\Omega}(e_1,e_3)=\overline{\Omega}(e_2,e_3)=0,  \nonumber
\end{eqnarray}
where
\begin{equation}\label{ab}
a(k+a^2+b^2)=ab=0 .
\end{equation}
If $a\neq 0$, from equations \eqref{3tripoles}-\eqref{ab} with
$\frak{a} = Span (f_1^3, f_2^3),$ one can show  that the Cartan
triple is smaller than $(\frak{o}(3), 0, -a^2\Omega_1),$ thus a
geometric realization has constant negative sectional curvature.
Note that in the case of a symmetric space $\bf M$ the
Ambrose-Singer connection and the Riemannian connection coincide as
connections on $\bf M,$ therefore one can read both the isometry
group and the curvature operator from
the Cartan triple ( see equations (2)).\\
If $a=0,$ a straightforward computation shows that a geometric
realization has the Ricci quadratic form given by
\begin{equation}\label{ricci}
Ric(x)=(k+b^2)((x^1)^2+(x^2)^2)+2b^2(x^3)^2.
\end{equation}
Since the sectional curvature is positive if and only the Ricci
quadratic form is positively defined, it follows that a geometric
realization of such a Cartan triple has positive sectional curvature
only if $b\neq 0$ and $k>-b^2.$ Moreover, Cartan \cite{C} showed
that in this case a geometric realization is topologically a three
dimensional sphere. Also note that such a Cartan triple is maximal
iff  $b\neq 0, b^2 \neq k>-b^2.$ From equation \eqref{3tripoles} it
follows that the isometry group of a geometric realization
 of such a Cartan triple is a $4$ dimensional extension
  of the the group  $\mathbb{S}^3$ of unit quaternions
 by the group of unit complex numbers $\mathbb{S}^1$. Indeed from (2) and (12),
with $a=0,$ it follows that
\begin{eqnarray}\label{3trip-4}
\Big[e_1,e_2 \Big]=-2be_3-(k+b^2)f_1^2 \nonumber \\
\Big[e_1,e_3 \Big]= be_2,
\Big[e_2,e_3 \Big]=-be_1, \\
\Big[f_1^2, e_1 \Big] = - e_2, \Big[f_1^2, e_2 \Big] = e_1,
\Big[f_1^2, e_3 \Big] = 0, \nonumber
\end{eqnarray}
showing that if $K$ is the simply connected Lie group of Lie algebra
$\frak{k}=:\frak{k}(\frak{g}_1, \Gamma,  \overline{\Omega}),$ then
$\frak{s}=:Span(e_1, e_2, 2be_3+(k+b^2)f_1^2)$ is a Lie subalgebra
of $\frak{k}$ that is transverse to $\frak{g}_1.$ From the general
Cartan triple approach \cite{P1} it follows that $\frak{s}$ is
isomorphic to a {\it transitive Killing algebra} of the geometric
realization $\bf M$ of $(\frak{g}_1, \Gamma,  \overline{\Omega}), $
that is $\frak{s}$ can be identified with a Lie algebra of Killing
vector fields on $\bf M,$ such that for any point $x \in M,$ the
evaluation map $ev_x : \frak{s} \rightarrow T_x\bf M, ev_x(\xi) =
\xi(x)$ is onto. Since $dim\frak{s}=dim(\bf M),$ it follows that
$ev_x$ is an isomorphism. Thus the connected Lie subgroup $S$ of
$\bf K$ of Lie subalgebra $\frak{s}$ acts transitively on $M,$ with
a discrete isotropy group (fiber $F$ of the projection $\pi:S
\rightarrow M, \pi(g)=g(x)$). Note that from the structural
equations of this Lie subalgebra, it turns out that $S$ is
isomorphic to $\mathbb{S}^3,$ the group of unit quaternions. From
the exact homotopy sequence of this fibration, since $S$ is
connected and simply connected we get $ \pi_0(F) =0,$ that is $F$ is
the fiber is the trivial subgroup of $S$, showing that $\pi$ is a
diffeomorphism. Note that  and since $\bf K$ is compact, the
isotropy group (in this case the Lie subgroup of $\bf K$ whose Lie
algebra is $\mathbb{R}f_1^2$) is a closed one dimensional
subgroup of a compact Lie group and therefore it is isomorphic to $\mathbb{S}^1.$  \\
\indent \'Elie Cartan (\cite{C}, p.305) integrated the structural
equations (15) and calculated the Riemannian structure of the
corresponding homogeneous spaces.

\begin{definition}
A geometric realization of $(\frak{g}_1, \Gamma,
\overline{\Omega}),$ with  $a=0, b\neq 0$ and $b^2 \neq k>-b^2$ is
said to be a Cartan sphere.
\end{definition}
Using (15) one can construct an isomorphism of Lie algebras \\
$\phi: \frak{k}(\frak{g}_1, \Gamma,  \overline{\Omega}) \rightarrow
\frak{u}(2),$ with $\phi(\frak{s})=\frak{su}(2).$ Shankar \cite{S}
mentions that $\mathbb{S}^3$ (regarded here as $SU(2)$ ) admits
Riemannian homogeneous structures ( not normal ) with $U(2)$ as the
connected component of the full group of isometries. These are the
only possible $4$ dimensional full groups of isometries of simply
connected three dimensional  Riemannian  homogeneous manifolds with
positive sectional curvature. Therefore we have following result.
\begin{proposition}
The Cartan spheres correspond to the three dimensional Riemannian
homogeneous manifolds Shankar's list, whose isometry groups are
locally isomorphic to $U(2).$
\end{proposition}
At the end of this section Cartan obtained the following result.
\begin{theorem}\label{3dsph}(Cartan)
Any geometric realization of a Cartan triple in 3D that admits a
four dimensional transitive group of isometry is homeomorphic to the
Euclidean space, a product of a round sphere and an Euclidean line
or a 3D sphere. Those geometric realizations of positive scalar
curvature are spheres, but there are 3D spheres, geometric
realizations of 3D Cartan triples that do not have a positive scalar
curvature.
\end{theorem}
\begin{remark} Theorem \ref{3dsph} was later extended by
in \cite{P2} who showed that geometric realizations of 3D Cartan
triples, including triples with $\mathfrak g = 0,$ can be deformed
to one of the eight 3D geometries of Thurston, that are the building
blocks of 3D manifolds according to Thurston and Perelman.
\end{remark}
\section{ General intransitive groups of isometries. Groups whose orbits are curves or surfaces}
In the last two sections of Chapter XII, Cartan considers the case
of intransitive groups of isometries. If $\bf K,$ is an intransitive
group of isometries of an $n+k$ dimensional Riemannian manifold
$(\bf M, g),$ and $\bf K(x)$ be an orbit of maximum dimension $n$ of
$\bf K,$ then $\bf K$ is transitive on $\bf K(x),$ therefore Cartan
shows that one can select coordinates $x^1, \dots x^n, y^1, \dots,
y^r),$ locally around $\bf K(x)$ such that $y^a, a = 1,\dots, r$ are
invariants of the group $\bf K$ and for fixed $y^1, \dots, y^r),$
the group acts transitively on the submanifold described by $x^1,
\dots, x^n.$ This is a {\it foliation} of $\bf M$ by homogeneous
spaces. Finally in Chapter XII of \cite{C}, Cartan describes the
form of the metric tensor, in case the maximal orbits are curves or
surfaces.

\end{document}